\documentclass[12pt,a4paper,reqno]{amsart}

\usepackage{amssymb}
\usepackage{amscd}
\usepackage[pdftex,pdfpagelabels]{hyperref}
\usepackage{enumerate}
\usepackage{comment}
\usepackage{graphicx}
\usepackage{siunitx}
\usepackage{tikz-cd}
\usepackage{stix}
\usepackage{bm}
\DeclareMathAlphabet\mathbfcal{LS2}{stixcal}{b}{n}
\numberwithin{equation}{section}
\usepackage{cleveref}
\usepackage[pagewise]{lineno} % \linenumbers

\usepackage{mathtools}%                  http://www.ctan.org/pkg/mathtools
\usepackage[tableposition=top]{caption}% http://www.ctan.org/pkg/caption
\usepackage{booktabs,dcolumn}%           http://www.ctan.org/pkg/dcolumn + http://www.ctan.org/pkg/booktabs

\DeclareFontFamily{OT1}{rsfs}{}
\DeclareFontShape{OT1}{rsfs}{n}{it}{<-> rsfs10}{}
\DeclareMathAlphabet{\mathscr}{OT1}{rsfs}{n}{it}

\theoremstyle{plain}

\newtheorem{theorem}{Theorem}[section]
\newtheorem{proposition}[theorem]{Proposition}

\theoremstyle{definition}

\newtheorem{remark}[theorem]{Remark}

\renewcommand\P{\mathbb{P}}
\newcommand\E{\mathbb{E}}

\newcommand\Var{\mathrm{Var}}

\newcommand\F{\mathbf{F}}
\newcommand\FP{\mathbf{FP}}

\newcommand\Kakeya{{\operatorname{Kakeya}}}
\newcommand\Nikodym{{\operatorname{Nikodym}}}

\newcommand\eps{\varepsilon}

\begin{document}

\title{New Nikodym set constructions over finite fields}

\author{Terence Tao}
\address{UCLA Department of Mathematics, Los Angeles, CA 90095-1555.}
\email{tao@math.ucla.edu}


\begin{abstract}  For any fixed dimension $d \geq 3$ we construct a Nikodym set in $\F_q^d$ of cardinality $q^d - 
(\frac{d-2}{\log 2} +1+o(1)) q^{d-1} \log q$ in the limit $q \to \infty$, when $q$ is an odd prime power.  This improves upon the naive random construction, which gives a set of cardinality
$q^d - (d-1+o(1)) q^{d-1} \log q$, and is new in the regime where $\F_q$ has unbounded characteristic and $q$ not a perfect square.  While the final proofs are completely human generated, the initial ideas of the construction were inspired by output from the tools \texttt{AlphaEvolve} and \texttt{DeepThink}.  We also present a simple construction of Nikodym sets in $\F_q^2$ for $q$ a perfect square that is a special case of known unital-based constructions, and matches the existing bounds of $q^2 - q^{3/2} + O(q \log q)$, assuming that $q$ is not the square of a prime $p \equiv 3 \pmod{4}$.
\end{abstract}  

\maketitle

\section{Introduction}

Let $q$ be a prime power, let $d \geq 1$ a fixed dimension, and let $\F_q^d$ be the $d$-dimensional vector space over the finite field $\F_q$ of order $q$.  We will be interested here in the regime where $d$ is fixed and $q$ is large; in particular, the asymptotic notation we introduce in \Cref{notation-sec} will be adapted to this regime.  Define a \emph{direction} $\omega$ to be a set in $\F_q^d$ of the form
$$ \omega = [v_1,\dots,v_d] \coloneqq \{ (tv_1,\dots,tv_d): t \in \F_q \backslash \{0\} \}$$
for some $v_1,\dots,v_d \in \F_q^d$, not all zero, and define the projective space $\FP_q^{d-1}$ to be the space of all directions $\omega$, thus the cardinality $|\FP_q^{d-1}|$ of this space obeys the well-known formula
\begin{equation}\label{fpp-card}
|\FP_q^{d-1}| = \frac{q^d-1}{q-1} = q^{d-1} + O(q^{d-2}).
\end{equation}
Given a point $x \in \F_q^d$ and a direction $\omega = [v_1,\dots,v_d]$, the \emph{line} $\ell = \ell_{x,\omega}$ in $\F_q^d$ passing through $x$ in the direction $\omega$ is given by the formula
$$ \ell_{x,\omega} = x + (\omega \cup \{0\}) = \{ x + (tv_1,\dots,tv_d): t \in \F_q \}.$$
We also define the \emph{punctured line}
$$ \ell_{x,\omega}^* = x + \omega = \{ x + (tv_1,\dots,tv_d): t \in \F_q \backslash \{0\} \} = \ell_{x,\omega} \backslash \{x\}.$$
We recall two types of subsets of $\F_q^d$:
\begin{itemize}
    \item A \emph{Kakeya set} is a subset $K$ of $\F_q^d$ that contains a line in every direction, thus for all $\omega \in \F_q^{d-1}$ there exists $x \in \F_q^d$ such that $\ell_{x,\omega} \subset K$.
    \item A \emph{Nikodym set} is a subset $N$ of $\F_q^d$ that contains a punctured line through every point, thus for all $x \in \F_q^d$ there exists $\omega \in \FP_q^{d-1}$ such that $\ell_{x,\omega}^* \subset K$.
\end{itemize}
Let $\Kakeya(d,q)$ and $\Nikodym(d,q)$ denote the minimum cardinality of a Kakeya or Nikodym set in $\F_q^d$ respectively.  Determining the asymptotics of these quantities as $q \to \infty$ has been a topic of much attention in recent years.  By a standard projective transformation argument, one can relate these two quantities by the inequality
\begin{equation}\label{nikodym-lower}
    \mathrm{Nikodym}(d,q) \geq \mathrm{Kakeya}(d,q) - \frac{2q^{d-1}-q^{d-2}-q}{q^{d-1}-1} q^{d-1} \geq \mathrm{Kakeya}(d,q) - 2q^{d-1};
\end{equation}
for the convenience of the reader we give the simple derivation of this well-known argument\footnote{For the Euclidean analogue of this observation, see e.g., \cite[Theorem 11.11]{mattila} or \cite{tao}.}  in \Cref{proj-app}.  Thus, up to a lower order error, the smallest Nikodym set is at least as large as the smallest Kakeya set.

However, it is expected that there is a significant gap between the two quantities.  This is confirmed by the best known upper and lower bounds, which we briefly summarize here\footnote{Some of these references were originally located using the Gemini Deep Research tool, and subsequently verified by the author.}.
\begin{itemize}
    \item We trivially have $\Kakeya(1,q) = \Nikodym(1,q)=q$.
    \item $\Kakeya(2,q)$ is equal to $q(q+1)/2 + (q-1)/2$ when $q$ is odd and $q(q+1)/2$ when $q$ is even, so one has 
$$\Kakeya(2,q) = \frac{1}{2} q^2 + O(q)$$
in both cases \cite{block}.
    \item In contrast, from the theory of blocking sets, we have the lower bound
    $$ \Nikodym(2,q) \geq q^2 - q^{3/2} - 1 + \frac{1}{4}s(1-s)q,$$
    where $s$ is the fractional part of $\sqrt{q}$ \cite{szonyi}.  When $q$ is a perfect square, this bound is sharp up to a lower order error $O(q \log q)$ \cite{block-error}\footnote{In the notation of that paper, Nikodym sets are the ``green'' portion of a ``green--black coloring''.}, thus
\begin{equation}\label{perfect-square}
\Nikodym(2,q) = q^2 - q^{3/2} + O(q \log q)
\end{equation}
    in this case.  However, there is no obvious way to extend the upper bound component of \eqref{perfect-square} to the non-perfect-square case.
    \item In general, we have the bounds
    $$ \frac{q^d}{\left(2 - \frac{1}{q}\right)^{d-1}} \leq \Kakeya(d,q) \leq \frac{ q^d }{2^{d-1}}\left(1 + \frac{d+1-2^{-d+2}}{q} + O\left(\frac{1}{q^2}\right)\right);$$
    see \cite{bukh-chao}.  In particular, $\Kakeya(d,q) = \frac{ q^d }{2^{d-1}} + O(q^{d-1})$ and thus also
    $$\Nikodym(d,q) \geq \frac{q^d}{2^{d-1}}  + O(q^{d-1}),$$
    thanks to \eqref{nikodym-lower}.
    \item It is conjectured \cite[Conjecture 1.2]{lund} that
    \begin{equation}\label{nconj}
    \Nikodym(d,q) = q^d - o(q^d).
    \end{equation}
    By previous results, this is known for $d \leq 2$.     In the case of bounded characteristic (which in particular includes the case of even $q$) the stronger bound 
    \begin{equation}\label{qi}
    \Nikodym(d,q) = q^d - O(q^{(1-\eps)d})
    \end{equation}
    is known for some $\eps>0$ depending on $d$ and the characteristic \cite[Theorem 1.6]{guo}. In three dimensions, the conjecture \eqref{nconj} would be implied by a further conjecture on unions of lines \cite[Conjecture 1.4]{lund}.
    \item The classes of Kakeya and Nikodym sets can both be checked to be closed under Cartesian products, giving rise to the inequalities 
    $$\Kakeya(d_1+d_2,q) \leq \Kakeya(d_1,q) \Kakeya(d_2,q)$$
    and 
    $$\Nikodym(d_1+d_2,q) \leq \Nikodym(d_1,q) \Nikodym(d_2,q)$$ 
    for any $d_1,d_2 \geq 1$.  When $q$ is a perfect square, one can combine this observation with the constructions in \cite{block-error} (and the trivial bound $\Nikodym(1,q)=q$) to obtain an upper bound
    \begin{equation}\label{qtor}
    \Nikodym(d,q) \leq q^d - \left\lfloor \frac{d}{2} \right\rfloor q^{d-1/2} + O(q^{d-1} \log q)
    \end{equation}
    for any fixed $d \geq 1$.
\end{itemize}
We refer the reader to the cited papers for a longer discussion of the history of these problems and prior results.

The above results leave open the question of non-trivial upper bounds on $\Nikodym(d,q)$ for $d > 2$ (i.e., constructions of small Nikodym sets in three and higher dimensions) when $q$ has unbounded characteristic and is not a perfect square.  This question was explored in a recent collaboration \cite{gdm} using a number of modern tools, most notably the LLM-powered optimizer \texttt{AlphaEvolve} from Google Deepmind.  The outcomes of that exploration can be found at \href{https://github.com/google-deepmind/alphaevolve_repository_of_problems/tree/main/experiments/finite_field_nikodym_problem}{this repository} and can be summarized as follows.
\begin{itemize}
    \item By an application of \texttt{AlphaEvolve} (restricting to the case of prime $q$ for simplicity), an \href{https://github.com/google-deepmind/alphaevolve_repository_of_problems/blob/main/experiments/finite_field_nikodym_problem/finite_field_nikodym.ipynb}{algebraic construction was generated} for $d=3$ (in which one deleted a small number of low-degree algebraic varieties, and specifically $\{ (x,y,x^iy): x,y \in \F_q \backslash \{0\}$ for $0 < |i| \leq 4$) which numerically suggested the upper bound $\Nikodym(3,q) \leq q^3 - 8q^2$.
    \item An AI-generated proof of this bound \href{https://github.com/google-deepmind/alphaevolve_repository_of_problems/blob/main/experiments/finite_field_nikodym_problem/deep_think_docs/nikodym.tex}{was produced} by \texttt{DeepThink}, who generalized the construction (replacing $4$ with an arbitrary parameter to be optimized later), \href{https://github.com/google-deepmind/alphaevolve_repository_of_problems/blob/main/experiments/finite_field_nikodym_problem/deep_think_docs/nikodym2.tex}{eventually claiming} an upper bound of the form $\Nikodym(3,q) \leq q^3 - 2q^2 \log q + o(q^2 \log q)$.  However, as already acknowledged in the AI-generated output, the argument was only heuristic.
    \item A human inspection of the arguments then revealed that a simple application of the probabilistic method could be used to rigorously establish the bound 
    \begin{equation}\label{nik-easy}
    \Nikodym(d,q) \leq q^d - (d-1+o(1))q^{d-1} \log q.
    \end{equation}
    We reproduce this calculation in \Cref{prob-method}.  A \href{https://github.com/google-deepmind/alphaevolve_repository_of_problems/blob/main/experiments/finite_field_nikodym_problem/deep_think_docs/nikodym1.tex}{separate invocation} of \texttt{DeepThink} was also able to reconstruct a reasonably complete (informal) proof of this bound.
    \item A further inspection of the heuristics suggested that, by removing random quadratic varieties instead of the original varieties $\{ (x,y,x^iy): x,y \in \F_q \backslash \{0\}$, one could theoretically improve \eqref{nik-easy} (for odd $q$ at least) to 
    \begin{equation}\label{nik-conj}
    \Nikodym(d,q) \leq q^d - \left(\frac{d-1}{\log 2}+o(1)\right) q^{d-1} \log q;
    \end{equation}
    see \Cref{heuristic-method}.    However, this bound was only heuristic.
    \item After several failed attempts (by both the human authors and \texttt{DeepThink}) to make this heuristic precise, a weaker bound
    \begin{equation}\label{nik-conj-2}
    \Nikodym(d,q) \leq q^d - \left(\frac{d-2}{\log 2}+1+o(1)\right) q^{d-1} \log q 
    \end{equation}
    was established for odd $q$, which improved upon \eqref{nik-easy} when $d \geq 3$.  (For even $q$, the bound \eqref{qi} is superior, and for $q$ a perfect square, the bound \eqref{qtor} is superior.)  In particular, this shows that purely random sets are not the optimal Nikodym set construction in this regime. 
\end{itemize}

In \Cref{new-bound} we give a complete (and human-generated) proof of \eqref{nik-conj-2}:

\begin{theorem}[New upper bound on Nikodym sets]\label{nik-upper}  If $d \geq 3$ and $q$ is an odd prime power, then \eqref{nik-conj-2} holds in the asymptotic limit $q \to \infty$.
\end{theorem}

As described above, this theorem was not directly obtained by computer-assisted methods; however, the explorations, partial results, and heuristics generated by those methods were instrumental in allowing the author to discover the construction and verification methods needed to establish the result.  On the other hand, we were unable to obtain the improved bound \eqref{nik-conj} despite the (computer-assisted) heuristic argument suggesting it, and pose it instead as an open conjecture.

Repeating the above computer-assisted paradigm in two dimensions, we also obtained a construction of Nikodym sets in the plane, though as it turns out this construction is a special case of a known construction based on unitals\footnote{We thank Ferdinand Ihringer for this observation.}; see Remark \ref{notnew}.  More precisely, in \Cref{2D-sec} we show

\begin{theorem}[A $2D$ Nikodym set construction]\label{alternate}
Let $q$ be an odd prime power that is a perfect square, with $q$ not the square of a prime $p \equiv 3 \pmod{4}$.  Then there exists a Nikodym set in $\F_q^2$ of cardinality $q^2 - q^{3/2} + O(q \log q)$, thus recovering (the upper bound portion of) \eqref{perfect-square} in this case.
\end{theorem}

This construction was \href{https://github.com/google-deepmind/alphaevolve_repository_of_problems/blob/main/experiments/finite_field_nikodym_problem/finite_field_nikodym.ipynb}{initially found} by \texttt{AlphaEvolve} and \href{https://github.com/google-deepmind/alphaevolve_repository_of_problems/blob/main/experiments/finite_field_nikodym_problem/deep_think_docs/nikodym7.tex}{verified} by \texttt{DeepThink}, but the proof we give in \Cref{2D-sec} is human-written.  As we discuss in Remark \ref{notnew}, this construction can also be reconstructed from known results on unitals, but we give a self-contained verification of the construction here for the convenience of the reader.

\subsection{Notation}\label{notation-sec}

For the rest of the paper, $q$ is understood to be an odd prime power.
We write $X \ll Y$, $Y \gg X$, or $X = O(Y)$ to denote a bound of the form $|X| \leq CY$, where $C$ is a constant depending only on the dimension $d$.  We also write $X \asymp Y$ for $X \ll Y \ll X$.  We use $X = o(Y)$ to denote a bound of the form $|X| \leq c_{d,\eps}(q) Y$, where $c_{d,\eps}(q)$ can depend on the parameters $d,\eps,q$ and goes to zero as $q \to \infty$ for any fixed choice of $d,\eps$.

\subsection{Acknowledgments}

The author was supported by the James and Carol Collins Chair, the Mathematical Analysis \& Application Research Fund, and by NSF grants DMS-2347850, and is particularly grateful to recent donors to the Research Fund.  He particularly thanks his coauthors Bogdan Georgiev, Javier G\'omez-Serrano, and Adam Zsolt Wagner for the highly productive and enjoyable collaboration \cite{gdm}, and for making available the outputs of that collaboration for the purposes of writing the current paper.  We thank Will Sawin and an anonymous contributor to the author's blog for corrections, and Ferdinand Ihringer for pointing out the connection of Theorem \ref{alternate} with the results in \cite{baker}.  We are particularly indebted to Will Sawin for pointing out a simplification to the proof of Theorem \ref{nik-upper}.

\section{The probabilistic method}\label{prob-method}

In this section we establish \eqref{nik-easy}.  Fix $d \geq 2$.  It suffices to show that for any sufficiently small $\eps>0$, one can find a Nikodym set $N$ of cardinality at most
\begin{equation}\label{npd}
 |N| \leq q^d - (d-1+O(\eps))q^{d-1} \log q
 \end{equation}
if $q$ is large enough.  We remark that the arguments in this section do not require $q$ to be odd.

Assume $q$ to be large. We select $N$ completely at random, with each $x \in \F_q^d$ lying in $N$ with an independent probability of
$$ \P(x \in N) = 1 - \frac{d-1+\eps}{q} \log q.$$
In particular, the indicator variables $1_{x \in N}$ are iid with mean $1 - \frac{d-1+\eps}{q} \log q$ and variance $O(\log q/q)$.

We recall (a special case of) Bennett's inequality \cite{bennett}: if $X_1,\dots,X_n \in [0,1]$ are iid random variables of mean $\mu$ and variance $\sigma$ then
\begin{equation}\label{bennett-ineq}
\P( |X_1+\dots+X_n - n\mu| \geq t)\leq 2 \exp\left( - n\sigma^2 h\left(\frac{t}{n\sigma^2}\right)\right)
\end{equation}
for any $\lambda>0$, where $h(u) \coloneqq (1+u)\log(1+u)-u$ (so in particular $h(u) \gg u^2$ for $u =O(1)$).  Applying this inequality to $|N| = \sum_{x \in \F_q^d} 1_{x \in N}$ and $t = \eps q^{d-1} \log q$, we have
\begin{align*}
n \sigma^2 &\asymp q^{d-1} \log q \\
\frac{t}{n\sigma^2} &\asymp \eps
\end{align*}
and thus
$$
\P( |N| < q^d - (d-1+2\eps)q^{d-1} \log q ) \ll \exp( - c_{\eps} q^{d-1} \log q ) $$
for some $c_{\eps}>0$ independent of $q$.

Now we consider the probability that $N$ is a Nikodym set.  Each punctured line $\ell^*_{x,\omega}$ has a probability of
$$ \left(1 - \frac{d-1+\eps}{q} \log q\right)^{q-1} \asymp q^{1-d-\eps}$$
of lying in $N$.  These events are independent in $\omega$ for fixed $x$, thus the probability that none of the $\ell^*_{x,\omega}$, $\omega \in \FP_q^{d-1}$ lie in $N$ is at most
$$ \left( 1 - c'_d q^{1-d-\eps} \right)^{|\FP_q^{d-1}|} \ll \exp( - c''_d q^{\eps} )$$
for some constants $c'_d, c''_d > 0$.  By the union bound, the probability that $N$ fails to be a Nikodym set is then at most $q^d \exp( - c''_d q^{\eps} )$.  Thus, for $q$ large enough, the above construction produces a Nikodym set obeying the required bound \eqref{npd} with positive probability, giving \eqref{nik-easy}.

As discussed below, a more complicated variant of this method was \href{https://github.com/google-deepmind/alphaevolve_repository_of_problems/blob/main/experiments/finite_field_nikodym_problem/deep_think_docs/nikodym.tex}{first discovered} by \texttt{DeepThink}.  After suggesting a purely random construction, \texttt{DeepThink} was \href{https://github.com/google-deepmind/alphaevolve_repository_of_problems/blob/main/experiments/finite_field_nikodym_problem/deep_think_docs/nikodym1.tex}{also able to reconstruct} most of the details of the above argument. 

\subsection{Heuristic arguments}\label{heuristic-method}

In our experiments from \cite{gdm}, \texttt{AlphaEvolve} suggested constructions of Nikodym sets in which one deleted various low degree varieties from $\F_q^d$.  We then used \texttt{DeepThink} to performed a \href{https://github.com/google-deepmind/alphaevolve_repository_of_problems/blob/main/experiments/finite_field_nikodym_problem/deep_think_docs/nikodym.tex}{heuristic analysis} of these constructions, which we paraphrase here.  Let $V$ be a ``typical'' variety in $\F_q^d$ of degree $D=O(1)$.  Then the intersection of this variety with a typical line should (after a linear transformation) look like the roots in $\F_q$ of a ``typical'' degree $D$ polynomial.  Assuming that the characteristic of $\F_q$ exceeds $D$, the Chebotarev density theorem predicts that the probability that there are no such roots should equal the probability $\delta_D$ that a random permutation on $D$ elements is a derangement (has no fixed points).  As is well known, one has
$$ \delta_D = 1 - \frac{1}{1!} + \frac{1}{2!} - \dots + \frac{(-1)^D}{D!}.$$
Meanwhile, the Lang--Weil inequality \cite{lang-weil} suggests that $V$ should have cardinality about $q^{d-1}$.  Thus, if one were to remove $k$ varieties $V_1,\dots,V_k$ from $\F_q^d$ of degrees $D_1,\dots,D_k$, the resulting set $N$ should have cardinality approximately
\begin{equation}\label{npk}
 N \approx q^d - k q^{d-1}
 \end{equation}
and the probability that a given line lies in $N$ should be roughly $\delta_{D_1} \dots \delta_{D_k}$.  The probabilistic arguments just given then suggest that one has a good chance of being a Nikodym set provided that
\begin{equation}\label{deltad}
 \delta_{D_1} \dots \delta_{D_k} \ggg q^{1-d}.
\end{equation}

In the AI-generated analysis provided by \texttt{DeepThink}, it was noted that $\delta_D$ oscillated between $1/3$ and $1/2$, but eventually converged to $1/e$.  Using the latter value $1/e$, this suggested that one could take $k \approx \log (q^{d-1})$, which when inserted into \eqref{npk} recovers the bound \eqref{nik-easy}.  However, the AI-generated analysis conceded that this argument was highly heuristic, and the available error estimates in standard tools such as the Lang--Weil inequality or Chebotarev density theorem were inadequate to make this calculation rigorous.

A human inspection of the arguments then suggested that one should instead take $D_1=\dots=D_k=2$ (i.e., to only use quadratic varieties), both to simplify the rigorous analysis and to allow the $\delta_{D_i}$ to attain the maximum value of $1/2$.  The above heuristics then suggested that one could increase $k$ to approximately $\log(q^{d-1}) / \log 2$, thus predicting the bound \eqref{nik-conj}.  

An \href{https://github.com/google-deepmind/alphaevolve_repository_of_problems/blob/main/experiments/finite_field_nikodym_problem/deep_think_docs/nikodym2.tex}{initial attempt} to use \texttt{DeepThink} to reproduce these heuristics was unsuccessful, and identified a geometric flaw in the construction: as quadratic polynomials of one variable will usually have $0$ or $2$ roots, but only rarely have just one (repeated) root, any non-tangent line $\ell_{x,\omega}$ through a given point $x$ in a quadratic variety $V$ would encounter a second point in that variety, thus significantly weakening the possibility for creating Nikodym sets by removing quadratic varieties due to the need to restrict attention to tangent lines (with \texttt{DeepThink} proposing the $k=2$ case as the limit of the method).

To obtain the intermediate rigorous bound \eqref{nik-conj-2}, we made further modifications to the construction which did not originate from any AI tools.  Firstly, after removing some randomly selected quadratic varieties, we randomly added back some points to $N$ to counteract the specific geometric obstacle mentioned above.  We also allowed for a small number of $x$ to not have any punctured line $\ell^*_{x,\omega}$ contained in the set, as one could repair those failures ``manually'' by adding some punctured lines to the set.  However, even with these alterations, we were still faced with the significant technical difficulty that the error bounds in the Lang--Weil theorem were too weak (roughly speaking, they were $O(1/\sqrt{q})$ times the size of the main term, whereas we needed an accuracy that was roughly of the form $o(1/q)$).  While in principle the methods of \'etale cohomology could be used to get around this obstacle, the amount of algebraic geometry needed to implement this seemed formidable.

It is still in principle possible that one could obtain the conjecture \eqref{nik-conj} by a sufficient injection of algebraic geometry tools.  However, we found a simpler approach, based on the projective linear invariance of the problem, that could be blended with the probabilistic method of \Cref{prob-method} to achieve the intermediate bound \eqref{nik-conj-2} without needing advanced results from algebraic geometry, instead relying on retaining a large number of available directions per lines at key steps of the process to keep the error terms below $o(1/q)$.  We turn to this argument next.

\section{Proof of new upper bound}\label{new-bound}

We now begin the proof of \Cref{nik-upper}.  Fix $d \geq 3$, let $\eps>0$ be sufficiently small, and assume $q$ odd and sufficiently large compared to $d$ and $\eps$.  It will suffice to construct a Nikodym set $N$ of cardinality
\begin{equation}\label{N-bound}
|N| \leq q^d - \left(\frac{d-2}{\log 2}+1\right) q^{d-1} \log q + O( \eps q^{d-1} \log q).
\end{equation}

The construction of $N$, and verification of its properties, is rather lengthy, but we can summarize the strategy as follows.
\begin{enumerate}
    \item By the standard deletion method (which, in our context, would be more accurately called an insertion method), we can allow a small number of failures in the Nikodym property, reducing matters to constructing a suitable ``almost Nikodym set'' $N'$: see \Cref{main-2}.
    \item By adapting the probabilistic method from the previous section, it will suffice to construct a slightly larger set $N''$ that is a ``robust almost Nikodym set'' in the sense that most points $x \in\F_q^d$ have \emph{many} punctured lines $\ell^*_{x,\omega}$ that lie in $N''$; see \Cref{main-3}.
    \item We construct $N''$ by randomly deleting some quadratic varieties from $\F_q^d$, and then adding back a small random set $W$; see \Cref{construct-sec}.  The small random set can be used to handle those $x$ which lie on one of the deleted varieties; the main case is when $x$ avoids all of the varieties.
    \item After some routine changes of variable, it suffices to obtain (with very high probability) a lower bound on the size $|E_1 \cap \dots \cap E_k|$ of an intersection of certain random quadratic subsets $E_1,\dots,E_k$ of $\FP_q^{d-1}$; see \Cref{quad-control}.
    \item A direct computation of the mean and variance of $|E_1 \cap \dots \cap E_k|$, after removing some inconvenient constraints on the quadratic polynomials defining $E_1,\dots,E_k$, then gives the claim.
\end{enumerate}

\subsection{Step 1: Reduction to constructing an almost Nikodym set}

We give details of each of the steps in the above strategy.  We first observe that it will suffice to construct an ``almost Nikodym set'' in which the Nikodym condition fails for a relatively small proportion (about $o(1/q)$) of points $x$.  Namely, we reduce to showing

\begin{theorem}[Construction of almost Nikodym set]\label{main-2} There exists a set $N'$ of cardinality
\begin{equation}\label{N'-bound}
|N'| \leq q^d - \left(\frac{d-2}{\log 2}+1\right) q^{d-1} \log q + O( \eps q^{d-1} \log q)
\end{equation}
with the property that for all but $O(\eps q^{d-1})$ points $x \in \F_q^d$, there is a punctured line $\ell^*_{x,\omega}$ that is contained in $N'$.
\end{theorem}

We now explain why \Cref{main-2} implies \eqref{nik-conj-2}.  Let $N'$ be as in \Cref{main-2}, thus there is a set $E \subset \F_q^d$ with $|E| \ll \eps q^{d-1}$ such that every $x \in \F^d \backslash E$ has a punctured line $\ell^*_{x,\omega}$ contained in $N'$.  For each $x \in E$, the complement of $N'$ has cardinality $O(q^{d-1} \log q)$, so by the pigeonhole principle there is a punctured line $\ell^*_{x,\omega_x}$ which intersects the complement of $N'$ in only $O(\log q)$ points.  If one adjoins all these punctured lines $\ell^*_{x,\omega_x}$ to $N'$, one obtains a new set $N$ obeying the bound \eqref{N-bound} while also being a Nikodym set, as desired.

\subsection{Step 2: Reduction to constructing a robust almost Nikodym set}

The next step is to (partially) use the probabilistic construction from \Cref{prob-method} to perform a tradeoff between the size of the almost Nikodym set and the number of punctured lines through a typical point that is contained in this almost Nikodym set. More precisely, we now reduce to showing

\begin{theorem}[Construction of robust almost Nikodym set]\label{main-3}  There exists a set $N''$ of cardinality
\begin{equation}\label{N''-bound}
|N''| \leq q^d - \frac{d-2}{\log 2} q^{d-1} \log q + O( \eps q^{d-1} \log q)
\end{equation}
with the property that for all but $O(\eps q^{d-1})$ points $x \in \F_q^d$, there are at least $q^{1+\eps^2}$ punctured lines $\ell^*_{x,\omega}$ that are contained in $N''$.
\end{theorem}

We now explain why \Cref{main-3} implies \Cref{main-2}.  With $N''$ as in \Cref{main-2}, let $N'$ be a random subset of $N''$ in which each element of $N''$ lies in $N'$ with an independent probability of $1 - \frac{\log q}{q}$.  By the same application of Bennett's inequality \eqref{bennett-ineq} used in \Cref{prob-method}, we see that \eqref{N'-bound} holds with probability
$$ 1 - O(\exp( - c_{d,\eps} q^{d-1} \log q ))$$
for some $c_{d,\eps}>0$ independent of $q$.  Also, any punctured line $\ell^*_{x,\omega}$ that was already in $N''$, would have a probability of
$$ \left(1 - \frac{\log q}{q}\right)^{q-1} \asymp \frac{1}{q}$$
of lying in $N'$, with these events being independent as $\omega$ varies (keeping $x$ fixed).  From this independence, we easily see that if $x$ had at least $q^{1+\eps^2}$ punctured lines $\ell^*_{x,\omega}$ in $N''$, then with probability at least 
$$ \left( 1 - \frac{c}{q} \right)^{q^{1+\eps^2}} \ll \exp( - c'q^{\eps^2} )$$
for some absolute constants $c,c'>0$, at least one of these lines will also lie in $N'$.  By the union bound, we thus see that with probability $1 - O(q^d  \exp( - c'q^{\eps^2} ))$, all the $x$ which had $q^{1+\eps^2}$ punctured lines $\ell^*_{x,\omega}$ in $N''$, will have at least one of these punctured lines in $N'$.  For $q$ large enough, we then obtain a set $N'$ obeying the properties required for \Cref{main-2} with positive probability.

\subsection{Step 3: Construction of the robust almost Nikodym set}\label{construct-sec}

It remains to prove \Cref{main-3}. We construct the robust almost Nikodym set by first removing several random quadratic varieties, and then restoring a few random points to evade the previously noted geometric obstacle. More precisely, we construct $N''$ as follows.

\begin{itemize}
    \item Set $k \coloneqq \lfloor (1-\eps) \frac{d-2}{\log 2} \log q \rfloor$, so in particular $2^k \asymp q^{(1-\eps) (d-2)}$.
    \item Select $k$ (inhomogeneous) absolutely irreducible\footnote{That is to say, the quadratics have degree exactly two, and do not split as the product of two linear factors in the algebraic closure of $\F_q$.} quadratic polynomials $Q_1,\dots,Q_k \colon \F_q^d \to \F_q$ with coefficients in $\F_q$, chosen uniformly and independently at random.
    \item For each $i=1,\dots,k$, let $V_i = \{Q_i = 0\}$ be the zero locus of $Q_i$.  Let $W$ be a random subset of $\F_q^d$ with each element of $\F_q^n$ lying in $W$ with probability $\eps$ (independently of each other and of the $Q_1,\dots,Q_k$).  We then set
    $$ N'' \coloneqq \left(\F_q^d \backslash \bigcup_{i=1}^k V_i\right) \cup W.$$
\end{itemize}

The requirement that each $Q_i$ be absolutely irreducible is not onerous.  The total number of quadratic polynomials in $d$ variables is $q^{(d+1)(d+2)/2}$.  The portion that have degree strictly less than two is merely $q^{d+1}$.  The portion that are products of two linear forms in $\F_q$ is $O( q^{2d+2-1} )$ (the $-1$ is due to the ability to transfer a scalar from one factor to another).  Similarly, the portion that are the product of a linear form in $\F_{p^2}$, its conjugate, and a scalar is also $O(q^{2d+2-1})$.  Thus all but $O(q^{2d+2-1}/q^{(d+1)(d+2)/2}) = O(1/q^{d(d-1)/2})$ of the quadratic polynomials in $d$ variables are absolutely irreducible. This failure rate of $O(1/q^{d(d-1)/2})$ will be negligible in practice because $d \geq 3$.

We claim:
\begin{itemize}
    \item[(i)]  With probability $1-o(1)$, the set $N''$ obeys the bound \eqref{N''-bound}.
    \item[(ii)]  The expected number of $x \in \F_q^n$ for which $N''$ contains fewer than $q^{1+\eps^2}$ punctured lines $\ell^*_{x,\omega}$ is $O( \eps q^{d-1})$.
\end{itemize}
From Markov's inequality we then see that $N$ obeys the requirements of \Cref{main-3} with positive probability.

\subsection{Step 3a: Lower bounding the size of the construction}

It remains to verify (i) and (ii). We begin with (i).  Taking complements, it suffices to show that $\bigcup_{i=1}^k V_i \backslash W$ has cardinality
$(1-O(\eps)) k q^{d-1}$ with probability $1-o(1)$. By the law of large numbers and the independence hypotheses, it will suffice to show that $\bigcup_{i=1}^k V_i$ has cardinality $(1-O(\eps)) k q^{d-1}$ with probability $1-o(1)$.  

By the Lang--Weil bound \cite{lang-weil} and the irreducibility of $Q_i$ we have $|V_i| = q^{d-1} + O(q^{d-3/2})$ for all $i$.  
By the Bonferroni inequalities
$$\sum_{i=1}^k |V_i| - \sum_{1 \leq i < j \leq k} |V_i \cap V_j|\leq  \left|\bigcup_{i=1}^k V_i\right| \leq \sum_{i=1}^k |V_i|$$
it thus will suffice to show that with probability $1-o(1)$, one has
\begin{equation}\label{vij}
 |V_i \cap V_j| \ll q^{d-2}
 \end{equation}
for all $1 \leq i < j \leq k$ (noting that $k = O(\log q)$.  

The Lang--Weil (or DeMillo--Lipton--Schwartz--Zippel) bound gives \eqref{vij} as long as $Q_i$ and $Q_j$ are not scalar multiples of each other, which by direct counting we know to be the case with probability $1-O(1/q)$ (say) for any given $i,j$.  Thus by the union bound, with probability $1-o(1)$ one has \eqref{vij} for all $i,j$, giving the claim (i).

\subsection{Step 4: Initial reductions for verifying the robust almost Nikodym property}

Now we prove (ii).  By linearity of expectation, it suffices to show that for each $x \in \F_q^d$, the probability that $N''$ contains  fewer than $q^{1+\eps^2}$ punctured lines $\ell^*_{x,\omega}$ is $O( \eps/q)$.  Because the distribution of the random set $N''$ is stationary (translation-invariant), it suffices to do this for $x=0$.  Equivalently, it suffices to show that with probability $1-O(\eps/q)$, there are at least $q^{1+\eps^2}$ directions $\omega \in \FP_q^{d-1}$ that avoid $\bigcup_{i=1}^k V_i \backslash W$.

The geometric obstacle mentioned in \Cref{heuristic-method} is activated when one of the $Q_i$ vanishes at the origin, as this then means that most directions $\omega$ will also encounter an additional element of $V_i$ besides the origin.  The additional random set $W$ is inserted purely to eliminate this obstacle.  More precisely, we now reduce to showing

\begin{proposition}[Avoiding quadratic varieties]\label{quadratic}  Suppose that the absolutely irreducible quadratic polynomials $Q_1,\dots,Q_k$ are conditioned to be non-vanishing at the origin.  Then with probability $1-O(1/q^{1+(d-2)\eps})$, there are at least $q^{1+2\eps^2}$ directions $\omega \in \FP_q^{d-1}$ that avoid $\bigcup_{i=1}^k V_i$.  Similarly if one replaces $k$ by $k-1$.
\end{proposition}

Let us see why this proposition establishes (ii).  Clearly it already handles the case when all the $Q_i$ do not vanish at the origin.  Since each $Q_i$ has a $O(1/q)$ chance of vanishing at the origin, we see from independence and the union bound the event where two or more $Q_i$ vanish has probability $O(k^2/q^2) = O(\eps/q)$, which is acceptable.  It remains to handle the case when exactly one of the $Q_i$ vanishes at the origin. By symmetry and the union bound, we may assume without loss of generality that $Q_k$ vanishes at the origin, provided that we improve our probability bound slightly from $1-O(\eps/q)$ to to $1 - O(\eps/kq)$.

Henceforth we condition to the event $Q_k(0)=0$. By the proposition (with $k$ replaced by $k-1$), with probability $1-O(1/q^{1+(d-2)\eps}) = 1-O(\eps/kq)$ it is already the case that there are at least $q^{1+2\eps^2}$ elements of $\FP_q^{d-1}$ that avoid $\bigcup_{i=1}^{k-1} V_i$.  Now condition $Q_1,\dots,Q_{k-1}$ to be fixed with this property.  For any direction $\omega \in \FP_q^{d-1}$, the quadratic polynomial $Q_k$ (which already vanishes at the origin), either vanishes identically on $\omega$, or has at most one further zero.  The former case only occurs with probability $O(1/q^2)$ per line ($Q_k$ has to vanish to second order in the direction $\omega$, thus creating two additional vanishing conditions $D_\omega Q_k(0) = D_\omega^2 Q_k(0) = 0$ beyond the event $Q_k(0)=0$ that is already being conditioned to), so by Markov's inequality, with probability $1-O(1/q^2)$ this happens for fewer than half of the $q^{1+2\eps^2}$ elements of $\FP_q^{d-1}$ under consideration.  Condition $V_k$ to be fixed so that this occurs, so that $W$ is now the only source of randomness.  For each remaining problematic direction $\omega$, of which there are at least $q^{1+2\eps^2}/2$, there is a probability at least $\eps$ that $\omega \cap V_k$ is contained in $W$, with these events being independent in $\omega$.  By the Bennett inequality \eqref{bennett-ineq} and a routine calculation, this containment will hold for at least $q^{1+\eps^2}$ values of $\omega$ with probability $1-O(\exp(-c_\eps q^{2\eps^2})) = 1 - O(\eps/kq)$ for some $c_\eps > 0$, giving the claim.

It remains to establish \Cref{quadratic}, which contains assertions for both $k$ and $k-1$.  We just prove the claim for $k$, as the claim for $k-1$ is proven just by replacing $k$ by $k-1$ throughout.  By scaling, we can normalize each $Q_i$ to equal $1$ at the origin.  So now we can assume that each $Q_i$ takes the form
$$ Q_i(x_1,\dots,x_d) = \sum_{1 \leq a \leq b \leq d} q_{ab,i} x_a x_b + \sum_{a=1}^d l_{a,i} x_a + 1$$
where the coefficients $q_{ab,i}, l_{a,i}$ are independent and uniformly distributed among $\F_q$, then conditioned on the requirement of absolute irreducibility.  A given direction $\omega = [v_1,\dots,v_d]$ of $\FP_q^{d-1}$ will then avoid $V_i$ if the univariate quadratic
$$ t \mapsto \sum_{1 \leq a \leq b \leq d} q_{ab,i} v_a v_b t^2 + \sum_{a=1}^d l_{a,i} v_a t + 1$$
has no zeroes in $\F_q$, which by the quadratic formula is equivalent to the discriminant
$$ \tilde Q_i(v_1,\dots,v_d) := \left(\sum_{a=1}^d l_{a,i} v_a\right)^2 - 4\sum_{1 \leq a \leq b \leq d} q_{ab,i} v_a v_b  $$
being a quadratic non-residue.  (Note that the truth of this assertion is unchanged if we multiply $(v_1,\dots,v_d)$ by a non-zero scalar, so by abuse of notation we can also say here that $\tilde Q(\omega)$ is a quadratic non-residue.)

Observe that if $Q_i$ were reducible, thus $Q_i(x) = (1+A_i(x)) (1+B_i(x))$ for some homogeneous linear forms $A_i,B_i$ (with coefficients in $\F_{q^2}$), then the discriminant would be a perfect square:
$$ \tilde Q_i(v) = (A_i(v)+B_i(v))^2 - 4 A_i(v) B_i(v) = (A_i(v)-B_i(v))^2.$$
Conversely, if $\tilde Q_i$ were a perfect square $\tilde Q_i(v) = C_i(v)^2$ for some homogeneous linear form $C_i$ (again with coefficients in $\F_{q^2}$), then writing $\sum_{a=1}^d l_{a,i} x_a$ as $L_i(x)$, we have
$$ Q_i(x) = 1 + L_i(x) + \frac{1}{4} (L_i(x)^2 - C_i(x)^2) = \left(1 + \frac{1}{2} L_i(x) + C_i(x)\right) \left(1 + \frac{1}{2} L_i(x) - C_i(x)\right).$$
Thus the requirement that $Q_i$ be absolutely irreducible is equivalent to the requirement that $\tilde Q_i$ is not a perfect square (in the absolute sense).  Furthermore, since $Q_i$ can be recovered from $\tilde Q_i$ together with an arbitrary choice of coefficients $l_{a,i}$, we see that the uniform distribution of $Q_i$ amongst absolutely irreducible polynomials implies the uniform distribution of $\tilde Q_i$ amongst homogeneous quadratics that are not a perfect square.  Furthermore, the $\tilde Q_i$ are independent in $i$.   To summarize, we have reduced to showing

\begin{proposition}[Controlling intersection of quadratic sets]\label{quad-control}  Let $\tilde Q_1,\dots,\tilde Q_k$ be homogeneous quadratic polynomials on $d$ variables that are not perfect squares, chosen uniformly and independently at random amongst all such choices.  Then with probability $1 -O(1/q^{1+(d-2)\eps})$, one has
\begin{equation}\label{eik}
 |E_1 \cap \dots \cap E_k| \geq q^{1+2\eps^2},
 \end{equation}
where $E_i \subset \FP_q^{d-1}$ are the independent random sets
$$E_i := \{ \omega \in \FP_q^{d-1}: \tilde Q_i(\omega) \hbox{ is a non-residue}\}.$$
\end{proposition}

\subsection{Step 5: the second moment method}

We now give a simple proof of Proposition \ref{quad-control} that was \href{https://terrytao.wordpress.com/2025/11/12/new-nikodym-set-constructions-over-finite-fields/comment-page-1/#comment-689091}{provided to us} by Will Sawin, who also contributed Remark \ref{improv} below.  First, observe from direct counting that the proportion of homogeneous quadratics $\tilde Q$ that are perfect squares $L(x)^2$ is $O(1/q^{\frac{d(d+1)}{2}-d}) = O(1/q^{d(d-1)/2})$ (noting that the coefficients of $L$ must square to coefficients of $\tilde Q$).  For $d \geq 3$, this proportion is negligible for our purposes since $k/q^{d(d-1)/2} \ll 1/q^{1+(d-2)\eps}$.  Thus, to prove Proposition \ref{quad-control}, we may remove the requirement that $\tilde Q_1,\dots,\tilde Q_k$ not be perfect squares, so that the $\tilde Q_i$ are now uniformly distributed amongst all homogeneous quadratic polynomials.  Now we consider the random variable
$$ |E_1 \cap \dots \cap E_k| = \sum_{\omega \in \FP_q^{d-1}} \prod_{i=1}^k 1_{\tilde Q_i(\omega) \hbox{ is a non-residue}}.$$
For any $\omega = [v_1,\dots,v_d] \in \FP_q^{d-1}$, the random variable $\tilde Q_i(v_1,\dots,v_d)$ is uniformly distributed in $\F_q$, because $\tilde Q_i$ is now uniformly distributed amongst quadratic polynomials.  Hence the probability that $\tilde Q_i(\omega)$ is a quadratic non-residue is $1/2 + O(1/q)$.  By independence, the indicator random variable $\prod_{i=1}^k 1_{\tilde Q_i(\omega) \hbox{ is a non-residue}}$ thus has mean $(1/2+O(1/q))^k$.  For distinct $\omega = [v_1,\dots,v_d]$ and $\omega' =[v'_1,\dots,v'_d]$ in $\FP_q^{d-1}$, the random variables $\tilde Q_i(v_1,\dots,v_d)$ and $\tilde Q_i(v'_1,\dots,v'_d)$ are independent.  Thus the indicators $\prod_{i=1}^k 1_{\tilde Q_i(\omega) \hbox{ is a non-residue}}$ are pairwise independent.  This gives the first and second moment bounds
$$ \E |E_1 \cap \dots \cap E_k| = |\FP_q^{d-1}| \left(\frac{1}{2} + O\left(\frac{1}{q}\right)\right)^k = \left(1 + O\left(\frac{k}{q}\right)\right) 2^{-k} q^{d-1}$$
and
$$ \Var |E_1 \cap \dots \cap E_k| \leq |\FP_q^{d-1}| \left(\frac{1}{2} + O\left(\frac{1}{q}\right)\right)^k = \left(1 + O\left(\frac{k}{q}\right)\right) 2^{-k} q^{d-1}.$$
By choice of $k$, one has $(1+O(\frac{k}{q})) 2^{-k} q^{d-1} \asymp q^{1 + (d-2)\eps}$, and the claim now follows from Chebyshev's inequality.

\begin{remark}\label{improv} By calculating more moments beyond the second, it may be possible to obtain an analogue of Proposition \ref{quad-control} for larger values of $k$, and potentially recover the conjecture \eqref{nik-conj}.  The calculations become more subtle, however, as the $s$-wise independence between the events $\omega \in E_i$ eventually breaks down for $s$ large enough; we will not attempt to perform these calculations here.
\end{remark}

\section{A two-dimensional construction}\label{2D-sec}

In this section we prove \Cref{alternate}, which arose from a separate generation of \texttt{AlphaEvolve} and \texttt{DeepThink} in \cite{gdm}, though the constructions obtained are similar in that they both involve removing quadratic varieties from the entire space and then performing some probabilistic modifications.  The runs of \texttt{AlphaEvolve} in \cite{gdm} \href{https://github.com/google-deepmind/alphaevolve_repository_of_problems/blob/main/experiments/finite_field_nikodym_problem/finite_field_nikodym.ipynb}{generated multiple examples} of Nikodym sets formed by deleting a small number of parabolae of the form $\{ (x,y): y-x^2 = s\}$ for various parameters $s$, eventually settling on choosing $s$ to be three quarters of an ``imaginary axis'' $\sqrt{c} \F_{\sqrt{q}}$.  We \href{https://github.com/google-deepmind/alphaevolve_repository_of_problems/blob/main/experiments/finite_field_nikodym_problem/deep_think_docs/nikodym5.tex}{used} \texttt{DeepThink} to analyze this construction; this tool used the Weil bound and \href{https://github.com/google-deepmind/alphaevolve_repository_of_problems/blob/main/experiments/finite_field_nikodym_problem/deep_think_docs/nikodym8.tex}{optimized the parameters} to give a construction of shape $q^2 - q^{3/2} + O(q^{5/4} \log q)$.  In this paper, we replace the Weil bound analysis with a random construction to improve the error term from $O(q^{5/4} \log q)$ to $O(q \log q)$ to recover a (simple special case of) an existing construction, as we shall now describe.

We may assume that $q$ is large. We can view $\F_{\sqrt{q}}$ as a subfield of $\F_q$; since $\sqrt{q}$ is not a prime $p \equiv 3 \pmod{4}$, we see that $-1$ is a quadratic residue in $\F_{\sqrt{q}}$.  It is convenient to fix a quadratic non-residue $c$ of $\F_{\sqrt{q}}$ and a square root $\sqrt{c}$ in $\F_q$, thus every element in $\F_q$ can be uniquely represented as $a+b\sqrt{c}$ for $a,b \in \F_{\sqrt{q}}$, and we write $\mathrm{Re}(a+b\sqrt{c}) \coloneqq a$.  As $-1$ is a quadratic residue, $-c$ is a quadratic nonresidue.  Standard calculations (counting points on a hyperbola) then show that the quadratic form $a^2+cb^2$ for $a,b \in \F_{\sqrt{q}}$ attains the value $0$ exactly once (when $a=b=0$) and attains every other value in $\F_{\sqrt{q}}$ exactly $\sqrt{q}+1$ times. Writing $a^2+cb^2 = \mathrm{Re}( (a+b\sqrt{c})^2)$, we conclude that the expression $\mathrm{Re}(t^2)$ for $t \in \F_q$ attains the value $0$ exactly once (when $t=0$) and attains every other value in $\F_{\sqrt{q}}$ exactly $\sqrt{q}+1$ times. 

We consider the preliminary set
$$ N_0 \coloneqq \{ (x,y) \in \F_q^2: \mathrm{Re}(y-x^2) \neq 0 \}.$$
This is $\F_q^2$ with $\sqrt{q}$ parallel parabolas of the form $\{(x,y): y-x^2=s\}$ removed, so has cardinality $q^2 - q^{3/2}$.  We note that this set is nearly a Nikodym set in the following sense:
\begin{itemize}
\item[(i)] If a point $p$ lies outside $N_0$, then there is a punctured line $\ell^*_{p,\omega}$ through $p$ that is contained in $N_0$.
\item[(ii)] If instead $p$ lies in $N_0$, then are $\sqrt{q}+1$ punctured lines $\ell^*_{p,\omega}$ through $p$ that are in contained in $N_0$ except at one point.
\end{itemize}
To see these claims, write $p = (x_0,y_0)$ and consider directions of the form $\omega = [1,m]$ for some $m \in \F_{\sqrt{q}}$, then
$$ \ell^*_{p,\omega} = \{ (x_0+t, y_0+mt): t \in \F_p \backslash \{0\} \}$$
and a given point $(x_0+t, y_0+mt)$ on this punctured line will lie in $N_0$ unless
$$ \mathrm{Re}( y_0+mt - (x_0+t)^2 ) = 0$$
which we rearrange as
\begin{equation}\label{nat}
 A + \mathrm{Re}(m' t) - \mathrm{Re}(t^2) = 0
\end{equation}
where $A \coloneqq \mathrm{Re}(y_0-x_0^2)$ and $m' \coloneqq m - 2x_0$.

First suppose that $p$ does not lie in $N_0$, then $A=0$.  If we choose the slope $m$ to be $2x_0$, then $m'=0$, and so \eqref{nat} simplifies to $\mathrm{Re}(t^2)=0$.  As already discussed, this equation is only attained when $t=0$, giving (i).

Now suppose that $p$ lies in $N_0$, so $A \neq 0$, thus $A$ is of the form $A = -\mathrm{Re}(t_0^2)$ for $\sqrt{q}+1$ choices of $t_0 \in \F_q \backslash \{0\}$.  For each such $t_0$, we take $m = 2x_0 - 2t_0$, then \eqref{nat} can be rewritten as $\mathrm{Re}((t-t_0)^2)=0$, which then has a single solution at $t=t_0$, giving the claim (ii).

We then enlarge $N_0$ to a Nikodym set $N$ by a standard probabilistic construction, with each point $p$ not in $N_0$ lying in $N$ with an independent probability of $C \log q / \sqrt{q}$ for some constant $C>2$.  From Markov's inequality we see that $N$ has cardinality $q^2 - q^{3/2} + O(q \log q)$ with probability $\gg 1$.  On the other hand, by claim (i), every point $p$ outside of $N_0$ already has a punctured line $\ell^*_{p,\omega}$ in $N$, and by claim (ii), each point $p$ in $N_0$ will also have a punctured line
in $\ell^*_{p,\omega}$ with probability
$$ 1 - \left( 1 - C \log q / \sqrt{q} \right)^{-\sqrt{q}-1},$$
which is $1-o(1/q^2)$ since $C>2$.  By the union bound, we obtain a Nikodym set of the desired cardinality with positive probability.

\begin{remark}\label{notnew}  This construction is a special case of the probabilistic construction in \cite[\S 4.1]{block-error}, where the (non-classical) unital in question\footnote{While much of the paper \cite{block-error} is concerned with the classical Hermitian unitals, the probabilistic construction given in \cite[\S 4.1]{block-error} applies to arbitrary unitals, including the one under consideration here.}
is taken to be the complement of $N_0$ (the union of parallel parabolae), together with a point at infinity.  The fact that this is indeed a unital (which essentially amounts to verifying properties (i) and (ii) above) can be deduced from the construction in \cite{baker} (also reproduced in \cite[Result 1]{barwick}), after setting the $\beta$ parameter in \cite[Result 1]{barwick} to zero.  We thank Ferdinand Ihringer for these references.
\end{remark}

\appendix

\section{A projective transformation argument}\label{proj-app}

In this section we prove \eqref{nikodym-lower}.  Let $N$ be a Nikodym set; it will suffice to construct a Kakeya set $K$ of cardinality
\begin{equation}\label{kb}
 |K| \leq |N| + \frac{2q^{d-1}-q^{d-2}-q}{q^{d-1}-1} q^{d-1}.
 \end{equation}
By hypothesis, there is a map $x \mapsto \omega_x$ from $\F_q^d$ to $\FP_q^{d-1}$ such that $\ell^*_{x,\omega_x} \subset N$ for all $x \in \F_q^d$.  Given a randomly chosen hyperplane $\pi$ in $\F_q^d$, the probability that $\omega_x$ is parallel to $\pi$ for a given $x$ is
$$ \frac{|\FP_q^{d-2}|}{|\FP_q^{d-1}|} = \frac{q^{d-2}-1}{q^{d-1}-1}.$$
Thus, the expected number of $x \in \F_q^d$ with $\omega_x$ parallel to $\pi$ is $\frac{q^{d-2}-1}{q^{d-1}-1} q^d$.  By the probabilistic method (or pigeonhole principle), we may thus find a plane $\pi$ for which $\omega_x$ is parallel to $\pi$ for at most $\frac{q^{d-2}-1}{q^{d-1}-1} q^d$ choices of $x \in \F_q^d$.  Foliating $\F_q^d$ into $p$ translates of $\pi$ and applying the pigeonhole principle again, we may thus find one of these translates $\pi'$ for which at most $\frac{q^{d-2}-1}{q^{d-1}-1} q^{d-1}$ of the lines $\ell_{x,\omega_x}$ lie in $\pi'$.

By applying a general linear transformation (which does not affect the property of being a Nikodym set), we may assume $\pi'$ to be the hyperplane $\F_q^{d-1} \times \{0\}$.  Thus, there is an exceptional set $E \subset \F_q^{d-1}$ of cardinality
\begin{equation}\label{eb}
|E| \leq \frac{q^{d-2}-1}{q^{d-1}-1} q^{d-1}
\end{equation}
such that for all $x \in \F_q^{d-1}$, the direction $\omega_{(x,0)}$ is not horizontal, and thus of the form $[v_x,1]$ for some $v_x \in \F_q^{d-1}$.  From the definition of a Nikodym set, this means that
\begin{equation}\label{ntx}
 (x+tv_x, t) \in N 
\end{equation}
whenever $x \in \F_q^{d-1} \backslash E$ and $t \in \F_q \backslash \{0\}$.

Now we introduce the set
\begin{align*}
K &\coloneqq \left\{ \left(\frac{x}{t},\frac{1}{t}\right): (x,t) \in N, t \neq 0 \right\}\\
&\quad \cup \F_q^{d-1} \times \{0\} \\
&\quad \cup \{ (tx,t): x \in E, t \neq 0 \}.
\end{align*} 
From the union bound one has
$$ |K| \leq |N| + q^{d-1} + |E| (q-1)$$
which gives \eqref{kb} thanks to \eqref{eb} and a brief calculation.  Now we check that $K$ is a Kakeya set, thus we need to show it contains a line $\ell_{x_\omega,\omega}$ in every direction $\omega = [v_1,\dots,v_d] \in \FP_q^{d-1}$.  If $\omega$ is horizontal (i.e., $v_d=0$) we can simply take $x_\omega=0$ since $K$ contains $\F_q^{d-1} \times \{0\}$.  Similarly if $\omega$ is of the form $[x,1]$ for $x \in E$, since $K$ contains the origin as well as $\{ (tx,t): t \neq 0\}$.  The only remaining case is if $\omega = [x,1]$ for some $x \in \F_q^{d-1} \backslash E$.  But from \eqref{ntx} we see that $K$ contains $(\frac{x}{t} + v_x, \frac{1}{t})$ for all $t \in \F_q \backslash \{0\}$ as well as $(v_x,0)$, and so we can take $x_\omega \coloneqq (v_x,0)$ in this case.  This concludes the proof.

\end{document}